\documentclass[12pt,twoside]{article}

\date{} 

\usepackage{amssymb,amsmath,latexsym,amsthm,amsfonts}

\newcommand{\N}{\mbox{$I \kern -4pt N$}}      

\newcommand{\Z}{\mbox{$Z \kern -7.5pt Z$}}     

\newcommand{\Q}{\mbox{$Q \kern -8pt I$}}      

\newcommand{\R}{{\bf R}}
\newcommand{\C}{\mbox{$C \kern -8pt I$}}      

\newcommand{\vp}{\varphi}
\newcommand{\dxt}{dx\,dt}
\newcommand{\tell}{\tilde{\ell}}

\font\dsrom=dsrom10 scaled 1200
\def \1{\textrm{\dsrom{1}}}

\def\om       {\omega}

\def\Om       {\Omega}

\begin{document}

\centerline{} 

\centerline{} 

\centerline {\Large{\bf Local null controllability of the N-dimensional}} 

\centerline{} 

\centerline{\Large{\bf Navier-Stokes system with N-1 scalar controls}} 

\centerline{} 

\centerline{\Large{\bf in an arbitrary control domain}} 

\centerline{} 

\centerline{\bf {Nicol\'as Carre\~no}} 

\centerline{} 

\centerline{Universit\'e Pierre et Marie Curie-Paris 6} 

\centerline{UMR 7598 Laboratoire Jacques-Louis Lions, Paris, F-75005 France} 

\centerline{ncarreno@ann.jussieu.fr} 


\centerline{} 

\centerline{\bf {Sergio Guerrero}} 

\centerline{} 

\centerline{Universit\'e Pierre et Marie Curie-Paris 6} 

\centerline{UMR 7598 Laboratoire Jacques-Louis Lions, Paris, F-75005 France} 

\centerline{guerrero@ann.jussieu.fr} 


\newtheorem{Theorem}{\quad Theorem}[section] 

\newtheorem{Definition}[Theorem]{\quad Definition} 

\newtheorem{Corollary}[Theorem]{\quad Corollary} 

\newtheorem{Proposition}[Theorem]{\quad Proposition} 

\newtheorem{Lemma}[Theorem]{\quad Lemma} 

\newtheorem{Example}[Theorem]{\quad Example} 

\newtheorem{Remark}[Theorem]{\quad Remark}

\numberwithin{equation}{section}

\begin{abstract} 
 In this paper we deal with the local null controllability of the $N-$ dimensional Navier-Stokes system with internal controls having one vanishing component. The novelty of this work is that no condition is imposed on the control domain.
\end{abstract} 

{\bf Subject Classification:} 35Q30, 93C10, 93B05  \\

{\bf Keywords:} Navier-Stokes system, null controllability, Carleman inequalities

\section{Introduction}

Let $\Om$ be a nonempty bounded connected open subset of $\R^N$ ($N=2$ or $3$) of class $C^{\infty}$.
Let $T>0$ and let $\om\subset\Om$ be a (small) nonempty open subset which is the control domain. We will 
use the notation $Q=\Om\times(0,T)$ and $\Sigma=\partial\Om\times (0,T)$. 

We will be concerned with the following controlled Navier-Stokes system:
\begin{equation}\label{eq:NS}
\left\lbrace \begin{array}{ll}
    y_t - \Delta y + (y\cdot \nabla)y + \nabla p = v\1_{\om} & \mbox{ in }Q, \\
    \nabla\cdot y = 0 & \mbox{ in }Q, \\
    y = 0 & \mbox{ on }\Sigma, \\
    y(0) = y^0 & \mbox{ in }\Om,
 \end{array}\right.
\end{equation}
where $v$ stands for the control which acts over the set $\om$.

The main objective of this work is to obtain the local null controllability of system \eqref{eq:NS} by means of $N-1$ scalar controls, i.e., we will prove the existence of a number $\delta>0$ such that, for every $y^0\in X$ ($X$ is an appropriate Banach space) satisfying
$$\|y^0\|_X \leq \delta,$$
and every $i\in\{1,\dots,N\}$, we can find a control $v$ in $L^2(\om\times(0,T))^N$ with $v_i\equiv 0$ such that the corresponding solution to \eqref{eq:NS} satisfies
$$y(T)=0 \mbox{ in }\Om.$$

This result has been proved in \cite{E&S&O&P-N-1} when $\overline{\om}$ intersects the boundary of $\Om$. Here, we remove this geometric assumption and prove the null controllability result for any nonempty open set $\om\subset\Om$. A similar result was obtained in \cite{CorGue} for the Stokes system.

Let us recall the definition of some usual spaces in the context of incompressible fluids:
$$V=\{y\in H^1_0(\Om)^N:\,\nabla\cdot y=0 \mbox{ in }\Om  \}$$
and
$$H=\{y\in L^2(\Om)^N:\,\nabla\cdot y=0 \mbox{ in }\Om,\,y\cdot n =0 \mbox{ on }\partial\Om\}.$$

Our main result is given in the following theorem:

\begin{Theorem}\label{teo:nullcontrol}
Let $i\in\{1,\dots,N\}$. Then, for every $T>0$ and $\om\subset\Om$, there exists $\delta>0$ such that, for every $y^0\in V$ satisfying
$$\|y^0\|_V\leq \delta,$$
we can find a control $v\in L^2(\om\times(0,T))^N$, with $v_i\equiv 0$, and a corresponding solution $(y,p)$ to \eqref{eq:NS} such that
$$y(T)=0,$$
i.e., the nonlinear system \eqref{eq:NS} is locally null controllable by means of $N-1$ scalar controls for an arbitrary control domain.
\end{Theorem}

\begin{Remark}\label{rem1}
For the sake of simplicity, we have taken the initial condition in a more regular space than usual. However, following the same arguments as in \cite{E&S&O&P} and \cite{E&S&O&P-N-1}, we can get the same result by considering $y^0\in H$ for $N=2$ and $y^0\in H\cap L^4(\Om)^3$ for $N=3$.
\end{Remark}

To prove Theorem \ref{teo:nullcontrol}, we follow a standard approach (see for instance \cite{OlegN-S},\cite{E&S&O&P} and \cite{E&S&O&P-N-1}). We first deduce a null controllability result for a linear system associated to \eqref{eq:NS}:

\begin{equation}\label{eq:Stokes}
\left\lbrace \begin{array}{ll}
    y_t - \Delta y + \nabla p = f + v\1_{\om} & \mbox{ in }Q, \\
    \nabla\cdot y = 0 & \mbox{ in }Q, \\
    y = 0 & \mbox{ on }\Sigma, \\
    y(0) = y^0 & \mbox{ in }\Om,
 \end{array}\right.
\end{equation}
where $f$ will be taken to decrease exponentially to zero in $T$. We first prove a suitable Carleman estimate for the adjoint system of \eqref{eq:Stokes} (see \eqref{eq:adj-Stokes} below). This will provide existence (and uniqueness) to a variational problem, from which we define a solution $(y,p,v)$ to \eqref{eq:Stokes} such that $y(T)=0$ in $\Om$ and $v_i=0$. Moreover, this solution is such that $e^{C/(T-t)}(y,v)\in L^2(Q)^N\times L^2(\om\times (0,T))^N$ for some $C>0$.

Finally, by means of an inverse mapping theorem, we deduce the null controllability for the nonlinear system. 

This paper is organized as follows. In section 2, we establish all the technical results needed to deal with the controllability problems. In section 3, we deal with the null controllability of the linear system \eqref{eq:Stokes}. Finally, in section 4 we give the proof of Theorem \ref{teo:nullcontrol}.

\section{Some previous results}

In this section we will mainly prove a Carleman estimate for the adjoint system of \eqref{eq:Stokes}. In order to do so, we are going to introduce some weight functions. Let $\om_0$ be a nonempty open subset of $\R^N$ such that $\overline{\om_0}\subset \om$ and $\eta\in C^2(\overline{\Om})$ such that
\begin{equation}
 |\nabla \eta|>0 \mbox{ in }\overline{\Om}\setminus\om_0,\, \eta>0 \mbox{ in }\Om \mbox{ and } \eta \equiv 0 \mbox{ on }\partial\Om.
\end{equation}
The existence of such a function $\eta$ is given in \cite{FurIma}.
Let also $\ell\in C^{\infty}([0,T])$ be a positive function satisfying 
\begin{equation}
\begin{split}
&\ell(t) = t \quad \forall t \in [0,T/4],\,\ell(t) = T-t \quad \forall t \in [3T/4,T],\\
& \ell(t)\leq \ell(T/2),\, \forall t\in [0,T].
\end{split}
\end{equation}

Then, for all $\lambda\geq 1$ we consider the following weight functions:
\begin{equation}\label{pesos}
\begin{split}
&\alpha(x,t) = \dfrac{e^{2\lambda\|\eta\|_{\infty}}-e^{\lambda\eta(x)}}{\ell^{8}(t)},\, \xi(x,t)=\dfrac{e^{\lambda\eta(x)}}{\ell^8(t)},\\
&\alpha^*(t) = \max_{x\in\overline{\Om}} \alpha(x,t),\, \xi^*(t) = \min_{x\in\overline{\Om}} \xi(x,t),\\
&\widehat\alpha(t) = \min_{x\in\overline{\Om}} \alpha(x,t),\, \widehat\xi(t) = \max_{x\in\overline{\Om}} \xi(x,t).
\end{split}
\end{equation}
These exact weight functions were considered in \cite{ImaPuelYam}.

We consider now a backwards nonhomogeneous system associated to the Stokes equation:
\begin{equation}\label{eq:adj-Stokes}
\left\lbrace \begin{array}{ll}
    -\vp_t - \Delta \vp + \nabla \pi = g & \mbox{ in }Q, \\
    \nabla\cdot \vp = 0 & \mbox{ in }Q, \\
    \vp = 0 & \mbox{ on }\Sigma, \\
    \vp(T) = \vp^T & \mbox{ in }\Om,
 \end{array}\right.
\end{equation}
where $g\in L^2(Q)^N$ and $\vp^T \in H$. 
Our Carleman estimate is given in the following proposition.

\begin{Proposition}\label{prop:Carleman}
There exists a constant $\lambda_0$, such that for any $\lambda>\lambda_0$ there exist two constants $C(\lambda)>0$ and $s_0(\lambda)>0$ such that for any  $i\in\{1,\dots,N\}$, any $g\in L^2(Q)^N$ and any $\vp^T \in H$, the solution of \eqref{eq:adj-Stokes} satisfies
\begin{equation}\label{eq:Carleman}
\begin{split}
s^4\iint\limits_{Q}  e^{-5s\alpha^*}  (\xi^*)^4 |\vp|^2 dx\,dt  
&\leq  C \left( \iint\limits_{Q}  e^{-3s \alpha^*}|g|^2 dx\,dt \right. \\
&\left.+ s^7 \sum_{j=1, j\neq i}^{N} \int\limits_0^T\int\limits_{\om}e^{-2s\widehat\alpha - 3s\alpha^*}(\widehat\xi)^7|\vp_j|^2 dx\,dt \right)
\end{split}
\end{equation}
for every $s\geq s_0$.
\end{Proposition}

The proof of inequality \eqref{eq:Carleman} is based on the arguments in \cite{CorGue}, \cite{E&S&O&P} and a Carleman inequality for parabolic equations with non-homogeneous boundary conditions proved in \cite{ImaPuelYam}. In \cite{CorGue}, the authors take advantange of the fact that the laplacian of the pressure is zero, but this is not the case here. Some arrangements of equation \eqref{eq:adj-Stokes} have to be made in order to follow the same strategy. More details are given below.

Before giving the proof of Proposition \ref{prop:Carleman}, we present some technical results. We first present a Carleman inequality proved in \cite{ImaPuelYam} for parabolic equations with nonhomogeneous boundary conditions.
To this end, let us introduce the equation
\begin{equation}\label{eq:heatnonhom}
 u_t - \Delta u = f_0 + \sum_{j=1}^N \partial_j f_j \mbox{ in }Q,
\end{equation}
where $f_0,f_1,\dots,f_N\in L^2(Q)$. We have the following result.

\begin{Lemma}\label{teo:Cnonhom}
There exists a constant $\widehat{\lambda_0}$ only depending on $\Om$, $\om_0$, $\eta$ and $\ell$ such that for any $\lambda>\widehat{\lambda_0}$ there exist two constants $C(\lambda)>0$ and $\widehat{s}(\lambda)$, such that for every $s\geq \widehat{s}$ and every $u\in L^2(0,T;H^1(\Om))\cap H^1(0,T;H^{-1}(\Om))$ satisfying \eqref{eq:heatnonhom}, we have
\begin{multline}\label{eq:Cnonhom}
 \dfrac{1}{s}\iint\limits_Q e^{-2s\alpha} \dfrac{1}{\xi}|\nabla u|^2 \dxt + s\iint\limits_Q e^{-2s\alpha} \xi |u|^2 \dxt \\
\leq C\left( s^{-\frac{1}{2}} \|e^{-s\alpha}\xi^{-\frac{1}{4}}u\|^2_{H^{\frac{1}{4},\frac{1}{2}}(\Sigma)} + s^{-\frac{1}{2}} \|e^{-s\alpha}\xi^{-\frac{1}{8}}u\|^2_{L^2(\Sigma)} \right.\\
 + \frac{1}{s^2}\iint\limits_Q e^{-2s\alpha}\frac{|f_0|^2}{\xi^2}\dxt 
 +\sum_{j=1}^N\iint\limits_Q e^{-2s\alpha}|f_j|^2 \dxt \\
\left.+ s\int\limits_0^T\int\limits_{\om_0}e^{-2s\alpha}\xi|u|^2 \dxt \right).
\end{multline}
\end{Lemma}

Recall that
$$\|u\|_{H^{\frac{1}{4},\frac{1}{2}}(\Sigma)}=\left(\|u\|^2_{H^{1/4}(0,T;L^2(\partial\Om))} + \|u\|^2_{L^{2}(0,T;H^{1/2}(\partial\Om))} \right)^{1/2}.$$

The next technical result is a particular case of Lemma 3 in \cite{CorGue}.

\begin{Lemma}\label{lemma1}
There exists $C>0$ depending only on $\Om$, $\om_0$, $\eta$ and $\ell$ such that, for every $T>0$ and every $u\in L^2(0,T;H^1(\Om))$,
\begin{multline}\label{eq:lemma1}
 s^3\lambda^2\iint\limits_Q e^{-2s\alpha}\xi^3|u|^2\dxt \\
\leq C \left( s \iint\limits_Q e^{-2s\alpha}\xi|\nabla u|^2\dxt + s^3\lambda^2 \int\limits_0^T\int\limits_{\om_0} e^{-2s\alpha}\xi^3|u|^2\dxt \right),
\end{multline}
for every $\lambda\geq \widehat{\lambda_1}$ and every $s\geq C$.
\end{Lemma}

\begin{Remark}
In \cite{CorGue}, slightly different weight functions are used to prove Lemma \ref{lemma1}. Namely, the authors take $\ell(t)=t(T-t)$. However, this does not change the result since the important property is that $\ell$ goes to $0$ polynomially when $t$ tends to $0$ and $T$.
\end{Remark}

The next lemma can be readily deduced from the corresponding result for parabolic equations in \cite{FurIma}.

\begin{Lemma}\label{lemma2}
 Let $\zeta(x)=\exp(\lambda\eta(x))$ for $x\in\Om$. Then, there exists $C>0$ depending only on $\Om$, $\om_0$ and $\eta$ such that, for every $u\in H^1_0(\Om)$,
\begin{multline}\label{eq:lemma2}
\tau^6\lambda^8\int\limits_{\Om}e^{2\tau\zeta}\zeta^6|u|^2 dx +  \tau^4\lambda^6\int\limits_{\Om}e^{2\tau\zeta}\zeta^4|\nabla u|^2 dx \\
\leq C \left( \tau^3\lambda^4\int\limits_{\Om}e^{2\tau\zeta}\zeta^3|\Delta u|^2 dx + \tau^6\lambda^8\int\limits_{\om_0}e^{2\tau\zeta}\zeta^6| u|^2 dx \right),
\end{multline}
for every $\lambda\geq \widehat{\lambda_2}$ and every $\tau\geq C$.
\end{Lemma}

The final technical result concerns the regularity of the solutions to the Stokes system that can be found in \cite{Lady} (see also \cite{Temam}).

\begin{Lemma}
For every $T>0$ and every $f\in L^2(Q)^N$, there exists a unique solution $u\in L^2(0,T;H^2(\Om)^N)\cap H^1(0,T;H)$ to the Stokes system
\begin{equation*}
\left\lbrace \begin{array}{ll}
    u_t - \Delta u + \nabla p = f  & \mbox{ in }Q, \\
    \nabla\cdot u = 0 & \mbox{ in }Q, \\
    u = 0 & \mbox{ on }\Sigma, \\
    u(0) = 0 & \mbox{ in }\Om,
 \end{array}\right.
\end{equation*}
for some $p\in L^2(0,T;H^1(\Om))$, and there exists a constant $C>0$ depending only on $\Om$ such that
\begin{equation}\label{eq:regularity1}
\|u\|^2_{L^2(0,T;H^2(\Om)^N)} + \|u\|^2_{H^1(0,T;L^2(\Om)^N)}\leq C \| f\|^2_{L^2(Q)^N}.
\end{equation}
Furthermore, if $f\in L^2(0,T;H^2(\Om)^N)\cap H^1(0,T;L^2(\Om)^N)$, then\\ $u\in L^2(0,T;H^4(\Om)^N)\cap H^1(0,T;H^2(\Om)^N)$ and there exists a constant $C>0$ depending only on $\Om$ such that
\begin{equation}\label{eq:regularity2}
\begin{split}
\|u\|^2_{L^2(0,T;H^4(\Om)^N)} &+ \|u\|^2_{H^1(0,T;H^2(\Om)^N)}  \\ 
&\leq C( \| f\|^2_{L^2(0,T;H^2(\Om)^N)} + \| f\|^2_{H^1(0,T;L^2(\Om)^N)}).
\end{split}
\end{equation}
\end{Lemma}

\subsection{Proof of Proposition \ref{prop:Carleman}}

Without any lack of generality, we treat the case of $N=2$ and $i=2$. The arguments can be easily extended to the general case. We follow the ideas of \cite{CorGue}. In that paper, the arguments are based on the fact that $\Delta \pi=0$, which is not the case here (recall that $\pi$ appears in \eqref{eq:adj-Stokes}). For this reason, let us first introduce $(w,q)$ and $(z,r)$, the solutions of the following systems:
\begin{equation}\label{eq:u1}
\left\lbrace \begin{array}{ll}
    -w_t - \Delta w + \nabla q = \rho g & \mbox{ in }Q, \\
    \nabla\cdot w = 0 & \mbox{ in }Q, \\
    w = 0 & \mbox{ on }\Sigma, \\
    w(T) = 0 & \mbox{ in }\Om,
 \end{array}\right.
\end{equation}
and
\begin{equation}\label{eq:u2}
\left\lbrace \begin{array}{ll}
    -z_t - \Delta z + \nabla r = -\rho' \vp & \mbox{ in }Q, \\
    \nabla\cdot z = 0 & \mbox{ in }Q, \\
    z = 0 & \mbox{ on }\Sigma, \\
    z(T) = 0 & \mbox{ in }\Om,
 \end{array}\right.
\end{equation}
where $\rho(t)=e^{-\frac{3}{2}s\alpha^*}$. Adding \eqref{eq:u1} and \eqref{eq:u2}, we see that $(w+z,q+r)$ solves the same system as $(\rho \vp,\rho \pi)$, where $(\vp,\pi)$ is the solution to \eqref{eq:adj-Stokes}. By uniqueness of the Stokes system we have
\begin{equation}\label{u1+u2}
\rho\vp= w + z \mbox{ and }\rho\pi = q + r.
\end{equation}

For system \eqref{eq:u1} we will use the regularity estimate \eqref{eq:regularity1}, namely
\begin{equation}\label{eq:regularity}
\|w\|^2_{L^2(0,T;H^2(\Om)^2)} + \|w\|^2_{H^1(0,T;L^2(\Om)^2)}\leq C \|\rho g\|^2_{L^2(Q)^2},
\end{equation}
 and for system \eqref{eq:u2} we will use the ideas of \cite{CorGue}. Using the divergence free condition on the equation of \eqref{eq:u2}, we see that
$$\Delta r=0 \mbox{ in }Q.$$
Then, we apply the operator $\nabla\Delta=(\partial_1\Delta,\partial_2\Delta)$ to the equation satisfied by $z_1$ and we denote $\psi:=\nabla\Delta z_1$. We then have

$$-\psi_t - \Delta \psi= -\nabla(\Delta(\rho'\vp_1))\mbox{ in }Q.$$
We apply Lemma \ref{teo:Cnonhom} to this equation and we obtain
\begin{multline}\label{eq:Carl1}
 I(s;\psi):=\dfrac{1}{s}\iint\limits_Q e^{-2s\alpha} \dfrac{1}{\xi}|\nabla \psi|^2 \dxt + s\iint\limits_Q e^{-2s\alpha} \xi |\psi|^2 \dxt \\
\leq C\left( s^{-\frac{1}{2}} \|e^{-s\alpha}\xi^{-\frac{1}{4}}\psi\|^2_{H^{\frac{1}{4},\frac{1}{2}}(\Sigma)^2} + s^{-\frac{1}{2}} \|e^{-s\alpha}\xi^{-\frac{1}{8}}\psi\|^2_{L^2(\Sigma)^2} \right.\\
\left. +  \iint\limits_Q e^{-2s\alpha} |\rho'|^2|\Delta \vp_1|^2 \dxt + s\int\limits_0^T\int\limits_{\om_0}e^{-2s\alpha}\xi|\psi|^2 \dxt \right),
\end{multline}
for every $\lambda \geq \widehat{\lambda_0}$ and $s\geq \widehat{s}$.

We divide the rest of the proof in several steps:
\begin{itemize}
\item In Step 1, using Lemmas \ref{lemma1} and \ref{lemma2}, we estimate global integrals of $z_1$ and $z_2$ by the left-hand side of \eqref{eq:Carl1}.
\item In Step 2, we deal with the boundary terms in \eqref{eq:Carl1}.
\item In Step 3, we estimate all the local terms by a local term of $\vp_1$ and $\epsilon\, I(s;\vp)$ to conclude the proof.
\end{itemize}
Now, let us choose $\lambda_0 = \max\{\widehat{\lambda_0},\widehat{\lambda_1},\widehat{\lambda_2}\}$ so that Lemmas \ref{lemma1} and \ref{lemma2} can be applied and fix $\lambda\geq \lambda_0$. In the following, $C$ will denote a generic constant depending on $\Om$, $\om$ and $\lambda$.

\textbf{Step 1.} \underline{\textit{Estimate of $z_1$}}. We use Lemma \ref{lemma1} with $u=\Delta z_1$:
\begin{multline}\label{eq:step1-1}
 s^3\iint\limits_Q e^{-2s\alpha}\xi^3|\Delta z_1|^2\dxt \\
\leq C \left( s \iint\limits_Q e^{-2s\alpha}\xi|\psi|^2\dxt + s^3 \int\limits_0^T\int\limits_{\om_0} e^{-2s\alpha}\xi^3|\Delta z_1|^2\dxt \right),
\end{multline}
for every $s\geq C$.

Now, we apply Lemma \ref{lemma2} with $u=z_1\in H^1_0(\Om)$ and we get:
\begin{multline*}
\tau^6\int\limits_{\Om}e^{2\tau\zeta}\zeta^6|z_1|^2 dx +  \tau^4\int\limits_{\Om}e^{2\tau\zeta}\zeta^4|\nabla z_1|^2 dx \\
\leq C \left( \tau^3\int\limits_{\Om}e^{2\tau\zeta}\zeta^3|\Delta z_1|^2 dx + \tau^6\int\limits_{\om_0}e^{2\tau\zeta}\zeta^6| z_1|^2 dx \right),
\end{multline*}
for every $\tau\geq C$. Now we take
$$\tau=\frac{s}{\ell^8(t)}$$
for $s$ large enough so we have $\tau\geq C$. This yields to
\begin{multline*}
s^6\int\limits_{\Om}e^{2s\xi}\xi^6|z_1|^2 dx +  s^4\int\limits_{\Om}e^{2s\xi}\xi^4|\nabla z_1|^2 dx \\
\leq C \left( s^3\int\limits_{\Om}e^{2s\xi}\xi^3|\Delta z_1|^2 dx + s^6\int\limits_{\om_0}e^{2s\xi}\xi^6| z_1|^2 dx \right),\,t\in(0,T),
\end{multline*}
for every $s\geq C$. We multiply this inequality by
$$\exp\left( -2s\frac{e^{2\lambda\|\eta\|_{\infty}}}{\ell^8(t)} \right),$$
and we integrate in $(0,T)$ to obtain
\begin{multline*}
s^6\iint\limits_{Q}e^{-2s\alpha}\xi^6|z_1|^2 \dxt +  s^4\iint\limits_{Q}e^{-2s\alpha}\xi^4|\nabla z_1|^2 \dxt \\
\leq C \left( s^3\iint\limits_{Q}e^{-2s\alpha}\xi^3|\Delta z_1|^2 \dxt + s^6\int\limits_0^T\int\limits_{\om_0}e^{-2s\alpha}\xi^6| z_1|^2 \dxt \right),
\end{multline*}
for every $s\geq C$. Combining this with \eqref{eq:step1-1} we get the following estimate for $z_1$:
\begin{multline}\label{eq:step1-2}
s^6\iint\limits_{Q}e^{-2s\alpha}\xi^6|z_1|^2 dxdt +  s^4\iint\limits_{Q}e^{-2s\alpha}\xi^4|\nabla z_1|^2 dxdt + s^3\iint\limits_Q e^{-2s\alpha}\xi^3|\Delta z_1|^2 dxdt \\
\leq C \left( s \iint\limits_Q e^{-2s\alpha}\xi|\psi|^2\dxt + s^3 \int\limits_0^T\int\limits_{\om_0} e^{-2s\alpha}\xi^3|\Delta z_1|^2\dxt \right.\\
\left. +  s^6\int\limits_0^T\int\limits_{\om_0}e^{-2s\alpha}\xi^6| z_1|^2 \dxt\right),
\end{multline}
for every $s\geq C$.

 \underline{\textit{Estimate of $z_2$}}. Now we will estimate a term in $z_2$ by the left-hand side of \eqref{eq:step1-2}. From the divergence free condition on $z$ we find
\begin{equation}\label{eq:step1-3}
\begin{split}
s^4\iint\limits_Q e^{-2s\alpha^*}(\xi^*)^4 |\partial_2 z_2|^2 \dxt 
&= s^4\iint\limits_Q e^{-2s\alpha^*}(\xi^*)^4 |\partial_1 z_1|^2 \dxt \\
&\leq s^4\iint\limits_Q e^{-2s\alpha}\xi^4 |\nabla z_1|^2 \dxt.
\end{split}
\end{equation}
Since $z_2|_{\partial\Om}=0$ and $\Om$ is bounded, we have that
$$\int\limits_{\Om} |z_2|^2 dx\leq C(\Om)\int\limits_{\Om} |\partial_2 z_2| dx,$$
and because $\alpha^*$ and $\xi^*$ do not depend on $x$, we also have
$$s^4\iint\limits_Q e^{-2s\alpha^*}(\xi^*)^4 | z_2|^2 \dxt\leq C(\Om) s^4\iint\limits_Q e^{-2s\alpha^*}(\xi^*)^4 |\partial_2 z_2|^2 \dxt.$$
Combining this with \eqref{eq:step1-3} we obtain
\begin{equation}\label{eq:step1-4}
s^4\iint\limits_Q e^{-2s\alpha^*}(\xi^*)^4 | z_2|^2 \dxt\leq C s^4\iint\limits_Q e^{-2s\alpha}\xi^4 |\nabla z_1|^2 \dxt.
\end{equation}

Now, observe that by \eqref{u1+u2}, \eqref{eq:regularity} and the fact that $s^2 e^{-2s\alpha} (\xi^*)^{9/4}$ is bounded we can estimate the third term in the right-hand side of \eqref{eq:Carl1}. Indeed,
\begin{multline*}
\iint\limits_Q e^{-2s\alpha} |\rho'|^2|\Delta \vp_1|^2 \dxt = \iint\limits_Q e^{-2s\alpha} |\rho'|^2 |\rho|^{-2} |\Delta (\rho\vp_1)|^2 \dxt \\
\leq C \left( s^2\iint\limits_Q e^{-2s\alpha} (\xi^*)^{9/4}|\Delta w_1| \dxt + s^2\iint\limits_Q e^{-2s\alpha} (\xi^*)^{9/4}|\Delta z_1| \dxt \right) \\
\leq C \left(  \|\rho g\|^2_{L^2(Q)^2} + s^2\iint\limits_Q e^{-2s\alpha} (\xi^*)^{3}|\Delta z_1| \dxt \right).
\end{multline*}
Putting together \eqref{eq:Carl1}, \eqref{eq:step1-2}, \eqref{eq:step1-4} and this last inequality we have for the moment
\begin{multline}\label{eq:endstep1}
s^6\iint\limits_{Q}e^{-2s\alpha}\xi^6|z_1|^2 dxdt +  s^4\iint\limits_{Q}e^{-2s\alpha^*}(\xi^*)^4| z_2|^2 dxdt + s^3\iint\limits_Q e^{-2s\alpha}\xi^3|\Delta z_1|^2 dxdt \\
+\dfrac{1}{s}\iint\limits_Q e^{-2s\alpha} \dfrac{1}{\xi}|\nabla \psi|^2 \dxt + s\iint\limits_Q e^{-2s\alpha} \xi |\psi|^2 \dxt \\
\leq C \left( s^{-\frac{1}{2}} \|e^{-s\alpha}\xi^{-\frac{1}{4}}\psi\|^2_{H^{\frac{1}{4},\frac{1}{2}}(\Sigma)^2} + s^{-\frac{1}{2}} \|e^{-s\alpha}\xi^{-\frac{1}{8}}\psi\|^2_{L^2(\Sigma)^2}  \right.\\
+  \|\rho g\|^2_{L^2(Q)^2}
 + s\int\limits_0^T\int\limits_{\om_0}e^{-2s\alpha}\xi|\psi|^2 \dxt \\
\left. + s^3 \int\limits_0^T\int\limits_{\om_0} e^{-2s\alpha}\xi^3|\Delta z_1|^2\dxt+  s^6\int\limits_0^T\int\limits_{\om_0}e^{-2s\alpha}\xi^6| z_1|^2 \dxt\right),
\end{multline}
for every $s\geq C$.

\textbf{Step 2.} In this step we deal with the boundary terms in \eqref{eq:endstep1}.

 First, we treat the second boundary term in \eqref{eq:endstep1}. Notice that, since $\alpha$ and $\xi$ coincide with $\alpha^*$ and $\xi^*$ respectively on $\Sigma$,
\begin{equation*}
\begin{split}
\|e^{-s\alpha^*}\psi\|^2_{L^2(\Sigma)^2} &\leq C \|s^{\frac{1}{2}}e^{-s\alpha^*}(\xi^*)^{\frac{1}{2}}\psi\|_{L^2(Q)^2}\|s^{-\frac{1}{2}}e^{-s\alpha^*}(\xi^*)^{-\frac{1}{2}}\nabla \psi\|_{L^2(Q)^2} \\
& \leq C\left( s\iint\limits_Q e^{-2s\alpha^*}\xi^*|\psi|^2\dxt + \frac{1}{s}\iint\limits_Q e^{-2s\alpha^*}\frac{1}{\xi^*}|\nabla \psi|^2 \dxt \right),
\end{split}
\end{equation*}
so $\|e^{-s\alpha^*}\psi\|^2_{L^2(\Sigma)^2}$ is bounded by the left-hand side of \eqref{eq:endstep1}. On the other hand,
$$s^{-\frac{1}{2}}\|e^{-s\alpha}\xi^{-\frac{1}{8}}\psi\|^2_{L^2(\Sigma)^2}\leq C s^{-\frac{1}{2}} \|e^{-s\alpha}\psi\|^2_{L^2(\Sigma)^2},$$
and we can absorb $s^{-\frac{1}{2}}\|e^{-s\alpha}\psi\|^2_{L^2(\Sigma)^2}$ by taking $s$ large enough.

 Now we treat the first boundary term in the right-hand side of \eqref{eq:endstep1}. We will use regularity estimates to prove that $z_1$ multiplied by a certain weight function is regular enough. First, let us observe that from \eqref{u1+u2} we readily have
\begin{multline*}
 s^4 \iint\limits_Q e^{-2s\alpha^*}(\xi^*)^4|\rho|^2|\vp|^2 \dxt \\
\leq 2 s^4 \iint\limits_Q e^{-2s\alpha^*}(\xi^*)^4|w|^2 \dxt + 2 s^4 \iint\limits_Q e^{-2s\alpha^*}(\xi^*)^4|z|^2 \dxt.
\end{multline*}
Using the regularity estimate \eqref{eq:regularity} for $w$ 
we have 
\begin{multline}\label{eq:step2-2}
 s^4 \iint\limits_Q e^{-2s\alpha^*}(\xi^*)^4|\rho|^2|\vp|^2 \dxt \\
\leq C \left( \|\rho g\|^2_{L^2(Q)^2} + s^4 \iint\limits_Q e^{-2s\alpha^*}(\xi^*)^4|z|^2 \dxt\right),
\end{multline}
thus the term $\|s^2e^{-s\alpha^*}(\xi^*)^2\rho\vp\|^2_{L^2(Q)^2}$ is bounded by the left-hand side of \eqref{eq:endstep1} and $\|\rho g\|^2_{L^2(Q)^2}$.

We define now
$$\widetilde{z}:=se^{-s\alpha^*}(\xi^*)^{7/8}z,\,\widetilde{r}:=se^{-s\alpha^*}(\xi^*)^{7/8}r.$$
From \eqref{eq:u2} we see that $(\widetilde{z},\widetilde{r})$ is the solution of the Stokes system:
\begin{equation*}
\left\lbrace \begin{array}{ll}
    -\widetilde{z}_t - \Delta \widetilde{z} + \nabla \widetilde{r} = -se^{-s\alpha^*}(\xi^*)^{7/8}\rho' \vp - (se^{-s\alpha^*}(\xi^*)^{7/8})_t z & \mbox{ in }Q, \\
    \nabla\cdot \widetilde{z} = 0 & \mbox{ in }Q, \\
    \widetilde{z} = 0 & \mbox{ on }\Sigma, \\
    \widetilde{z}(T) = 0 & \mbox{ in }\Om.
 \end{array}\right.
\end{equation*}
Taking into account that
$$|\alpha^*_t| \leq C (\xi^*)^{9/8},\,|\rho'|\leq C s\rho (\xi^*)^{9/8}$$
and the regularity estimate \eqref{eq:regularity1} we have
\begin{multline*}
\|\widetilde{z}\|^2_{L^2(0,T;H^2(\Om)^2)\cap H^1(0,T;L^2(\Om)^2)} \\
\leq C \left( \|s^2e^{-s\alpha^*}(\xi^*)^2\rho\vp\|^2_{L^2(Q)^2} + \|s^2e^{-s\alpha^*}(\xi^*)^2 z\|^2_{L^2(Q)^2} \right),
\end{multline*}
thus, from \eqref{eq:step2-2}, $\|se^{-s\alpha^*}(\xi^*)^{7/8}z\|^2_{L^2(0,T;H^2(\Om)^2)\cap H^1(0,T;L^2(\Om)^2)}$ is bounded by the left-hand side of \eqref{eq:endstep1} and $\|\rho g\|^2_{L^2(Q)^2}$. From \eqref{u1+u2}, \eqref{eq:regularity} and this last inequality we have that
\begin{multline*}
\|se^{-s\alpha^*}(\xi^*)^{7/8}\rho\vp\|^2_{L^2(0,T;H^2(\Om)^2)\cap H^1(0,T;L^2(\Om)^2)} \\
\leq C \left( \|\rho g\|^2_{L^2(Q)^2} + \|\widetilde{z}\|^2_{L^2(0,T;H^2(\Om)^2)\cap H^1(0,T;L^2(\Om)^2)} \right),
\end{multline*}
and thus $\|se^{-s\alpha^*}(\xi^*)^{7/8}\rho\vp\|^2_{L^2(0,T;H^2(\Om)^2)\cap H^1(0,T;L^2(\Om)^2)}$ is bounded by the left-hand side of \eqref{eq:endstep1} and $\|\rho g\|^2_{L^2(Q)^2}$.

Next, let
$$\widehat{z}:=e^{-s\alpha^*}(\xi^*)^{-1/4}z,\,\widehat{r}:=e^{-s\alpha^*}(\xi^*)^{-1/4}r.$$
From \eqref{eq:u2}, $(\widehat{z},\widehat{r})$ is the solution of the Stokes system:
\begin{equation*}
\left\lbrace \begin{array}{ll}
    -\widehat{z}_t - \Delta \widehat{z} + \nabla \widehat{r} = -e^{-s\alpha^*}(\xi^*)^{-1/4}\rho' \vp - (e^{-s\alpha^*}(\xi^*)^{-1/4})_t z & \mbox{ in }Q, \\
    \nabla\cdot \widehat{z} = 0 & \mbox{ in }Q, \\
    \widehat{z} = 0 & \mbox{ on }\Sigma, \\
    \widehat{z}(T) = 0 & \mbox{ in }\Om.
 \end{array}\right.
\end{equation*}
From the previous estimates, it is not difficult to see that the right-hand side of this system is in $L^2(0,T;H^2(\Om)^2)\cap H^1(0,T;L^2(\Om)^2)$, and thus, using the regularity estimate \eqref{eq:regularity2}, we have
\begin{multline*}
\|\widehat{z}\|^2_{L^2(0,T;H^4(\Om)^2)\cap H^1(0,T;H^2(\Om)^2)} \\
\leq C \left( \|se^{-s\alpha^*}(\xi^*)^{7/8}\rho\vp\|^2_{L^2(0,T;H^2(\Om)^2)\cap H^1(0,T;L^2(\Om)^2)} \right. \\
\left.+ \|se^{-s\alpha^*}(\xi^*)^{7/8} z\|^2_{L^2(0,T;H^2(\Om)^2)\cap H^1(0,T;L^2(\Om)^2)} \right).
\end{multline*}
In particular, $e^{-s\alpha^*}(\xi^*)^{-1/4}\psi \in L^2(0,T;H^1(\Om)^2)\cap H^1(0,T;H^{-1}(\Om)^2)$  (recall that $\psi=\nabla\Delta z_1$) and
\begin{equation}\label{eq:step2-1}
\| e^{-s\alpha^*}(\xi^*)^{-1/4}\psi \|^2_{L^2(0,T;H^1(\Om)^2)} \mbox{ and } \| e^{-s\alpha^*}(\xi^*)^{-1/4}\psi \|^2_{H^1(0,T;H^{-1}(\Om)^2)}
\end{equation}
are bounded by the left-hand side of \eqref{eq:endstep1} and $\|\rho g\|^2_{L^2(Q)^2}$. 

To end this step, we use the following trace inequality
\begin{equation*}
\begin{split}
 &s^{-1/2}\|e^{-s\alpha}\xi^{-\frac{1}{4}}\psi\|^2_{H^{\frac{1}{4},\frac{1}{2}}(\Sigma)^2} = s^{-1/2}\|e^{-s\alpha^*}(\xi^*)^{-\frac{1}{4}}\psi\|^2_{H^{\frac{1}{4},\frac{1}{2}}(\Sigma)^2} \\
& \leq  C\,s^{-1/2} \left( \| e^{-s\alpha^*}(\xi^*)^{-1/4}\psi \|^2_{L^2(0,T;H^1(\Om)^2)} 
 + \| e^{-s\alpha^*}(\xi^*)^{-1/4}\psi \|^2_{H^1(0,T;H^{-1}(\Om)^2)} \right).
\end{split}
\end{equation*}
 By taking $s$ large enough in \eqref{eq:endstep1}, the boundary term $s^{-1/2}\|e^{-s\alpha}\xi^{-\frac{1}{4}}\psi\|^2_{H^{\frac{1}{4},\frac{1}{2}}(\Sigma)^2}$ can be absorbed by the terms in \eqref{eq:step2-1} and step 2 is finished.

Thus, at this point we have
\begin{multline}\label{eq:endstep2-2}
s^4\iint\limits_{Q}e^{-2s\alpha^*}(\xi^*)^4 |\rho|^2 |\vp|^2 \dxt + s^3\iint\limits_Q e^{-2s\alpha}\xi^3|\Delta z_1|^2\dxt \\
+\dfrac{1}{s}\iint\limits_Q e^{-2s\alpha} \dfrac{1}{\xi}|\Delta^2 z_1|^2 \dxt + s\iint\limits_Q e^{-2s\alpha} \xi |\nabla\Delta z_1|^2 \dxt \\
\leq C \left( \|\rho g\|^2_{L^2(Q)^2} 
+  s^6\int\limits_0^T\int\limits_{\om_0}e^{-2s\alpha}\xi^6| z_1|^2 \dxt \right. \\
\left. + s\int\limits_0^T\int\limits_{\om_0}e^{-2s\alpha}\xi|\nabla\Delta z_1|^2 \dxt
 + s^3 \int\limits_0^T\int\limits_{\om_0} e^{-2s\alpha}\xi^3|\Delta z_1|^2\dxt \right), \\
\end{multline}
for every $s\geq C$.

\textbf{Step 3.} In this step we estimate the two last local terms in the right-hand side of \eqref{eq:endstep2-2} in terms of local terms of $z_1$ and the left-hand side of \eqref{eq:endstep2-2} multiplied by small constants. Finally, we make the final arrangements to obtain \eqref{eq:Carleman}.

We start with the term $\nabla\Delta z_1$ and we follow a standard approach. Let $\om_1$ be an open subset such that $\om_0\Subset\om_1\Subset\om$ and let $\rho_1 \in C^2_c(\om_1)$ with $\rho_1\equiv 1$ in $\om_0$ and $\rho_1\geq 0$. Then, by integrating by parts we get
\begin{equation*}
\begin{split}
&s\int\limits_0^T\int\limits_{\om_0}e^{-2s\alpha}\xi|\nabla\Delta z_1|^2 \dxt \leq s\int\limits_0^T\int\limits_{\om_1}\rho_1 e^{-2s\alpha}\xi|\nabla\Delta z_1|^2 \dxt\\
&= -s\int\limits_0^T\int\limits_{\om_1}\rho_1 e^{-2s\alpha}\xi \Delta^2 z_1 \Delta z_1  \dxt + \frac{s}{2} \int\limits_0^T\int\limits_{\om_1} \Delta(\rho_1 e^{-2s\alpha}\xi)  |\Delta z_1|^2 \dxt.
\end{split}
\end{equation*}
Using Cauchy-Schwarz's inequality for the first term and
$$|\Delta(\rho_1 e^{-2s\alpha}\xi)|\leq C s^2e^{-2s\alpha}\xi^3,\, s\geq C$$
for the second one, we obtain for every $\epsilon > 0$
\begin{equation*}
\begin{split}
&s\int\limits_0^T\int\limits_{\om_0}e^{-2s\alpha}\xi|\nabla\Delta z_1|^2 \dxt \\
&\leq  \frac{\epsilon}{s} \int\limits_0^T\int\limits_{\om_1} e^{-2s\alpha}\frac{1}{\xi} |\Delta^2 z_1|^2  \dxt + C(\epsilon)s^3 \int\limits_0^T\int\limits_{\om_1} e^{-2s\alpha}\xi^3  |\Delta z_1|^2 \dxt,
\end{split}
\end{equation*}
for every $s \geq C$.

Let us now estimate $\Delta z_1$. Let $\rho_2 \in C^2_c(\om)$ with $\rho_2\equiv 1$ in $\om_1$ and $\rho_2\geq 0$. Then, by integrating by parts we get
\begin{multline*}
s^3\int\limits_0^T\int\limits_{\om_1}e^{-2s\alpha}\xi^3  |\Delta z_1|^2  \dxt 
\leq s^3\int\limits_0^T\int\limits_{\om}\rho_2 e^{-2s\alpha}\xi^3|\Delta z_1|^2 \dxt\\
= 2s^3\int\limits_0^T\int\limits_{\om} \nabla(\rho_2 e^{-2s\alpha}\xi^3) \nabla\Delta z_1 \cdot  z_1  \dxt
  +s^3\int\limits_0^T\int\limits_{\om}\Delta(\rho_2 e^{-2s\alpha}\xi^3)\Delta z_1 \cdot z_1 \dxt \\
  +s^3\int\limits_0^T\int\limits_{\om}\rho_2 e^{-2s\alpha}\xi^3 \Delta^2 z_1 \cdot z_1  \dxt.
\end{multline*}
Using
$$|\nabla(\rho_2 e^{-2s\alpha}\xi^3)|\leq C se^{-2s\alpha}\xi^4,\, s\geq C,$$
for the first term in the right-hand side of this last inequality,
$$|\Delta(\rho_2 s^3e^{-2s\alpha}\xi^3)|\leq C s^5e^{-2s\alpha}\xi^5,\, s\geq C,$$
for the second one and Cauchy-Schwarz's inequality we obtain for every $\epsilon > 0$
\begin{equation*}
\begin{split}
&s^3\int\limits_0^T\int\limits_{\om_1}e^{-2s\alpha}\xi^3|\Delta z_1|^2 \dxt \\
&\leq  \epsilon \,\left( \frac{1}{s} \int\limits_0^T\int\limits_{\om} e^{-2s\alpha} \frac{1}{\xi} | \Delta^2 z_1|^2  \dxt + s \int\limits_0^T\int\limits_{\om} e^{-2s\alpha}\xi |\nabla \Delta z_1|^2  \dxt \right. \\
&\left.+ s^3 \int\limits_0^T\int\limits_{\om} e^{-2s\alpha}\xi^3 |\Delta z_1|^2  \dxt\right) 
 + C(\epsilon)s^7 \int\limits_0^T\int\limits_{\om} e^{-2s\alpha}\xi^7  | z_1|^2 \dxt,
\end{split}
\end{equation*}
for every $s \geq C$.

Finally, from \eqref{u1+u2} and \eqref{eq:regularity} we readily obtain
\begin{equation*}
\begin{split}
& s^7 \int\limits_0^T\int\limits_{\om} e^{-2s\alpha}\xi^7  |z_1|^2 \dxt \\
& \leq 2 s^7 \int\limits_0^T\int\limits_{\om} e^{-2s\alpha}\xi^7 |\rho|^2 |\vp_1|^2 \dxt + 2 s^7 \int\limits_0^T\int\limits_{\om} e^{-2s\alpha}\xi^7  |w_1|^2 \dxt \\
& \leq 2 s^7 \int\limits_0^T\int\limits_{\om} e^{-2s\alpha}\xi^7 |\rho|^2 |\vp_1|^2 \dxt + C \|\rho g\|^2_{L^2(Q)^2}.
\end{split} 
\end{equation*}
This concludes the proof of Proposition \ref{prop:Carleman}.

\section{Null controllability of the linear system}

Here we are concerned with the null controllability of the system
\begin{equation}\label{eq:Stokes2}
\left\lbrace \begin{array}{ll}
    y_t - \Delta y + \nabla p = f + v\1_{\om} & \mbox{ in }Q, \\
    \nabla\cdot y = 0 & \mbox{ in }Q, \\
    y = 0 & \mbox{ on }\Sigma, \\
    y(0) = y^0 & \mbox{ in }\Om,
 \end{array}\right.
\end{equation}
where $y^0\in V$, $f$ is in an appropiate weighted space and the control 
$v\in L^2(\om\times (0,T))^N$ is such that $v_i=0$ for some $i\in \{1,\dots,N\}$.

Before dealing with the null controllability of \eqref{eq:Stokes2}, we will deduce a new Carleman inequality with weights not vanishing at $t=0$. To this end, let us introduce the following weight functions:
\begin{equation}\label{pesos2}
\begin{split}
&\beta(x,t) = \dfrac{e^{2\lambda\|\eta\|_{\infty}}-e^{\lambda\eta(x)}}{\tell^{8}(t)},\, \gamma(x,t)=\dfrac{e^{\lambda\eta(x)}}{\tell^8(t)},\\
&\beta^*(t) = \max_{x\in\overline{\Om}} \beta(x,t),\, \gamma^*(t) = \min_{x\in\overline{\Om}} \gamma(x,t),\\
&\widehat{\beta}(t) = \min_{x\in\overline{\Om}} \beta(x,t),\, \widehat{\gamma}(t) = \max_{x\in\overline{\Om}} \gamma(x,t),
\end{split}
\end{equation}
where
\begin{equation*}
\tell(t)= \left\lbrace
\begin{array}{ll}
\|\ell\|_{\infty} & 0\leq t \leq T/2, \\
\ell(t) & T/2< t \leq T. 
\end{array}
\right.
\end{equation*}

\begin{Lemma}\label{lemma:Carleman2}
Let $i\in\{1,\dots,N\}$ and let $s$ and $\lambda$ be like in Proposition \ref{prop:Carleman}. Then, there exists a constant $C>0$ (depending on $s$ and $\lambda$) such that every solution $\vp$ of \eqref{eq:adj-Stokes} satisfies:
\begin{multline}\label{eq:Carleman2}
\iint\limits_Q e^{-5s\beta^*}(\gamma^*)^4|\vp|^2 \dxt + \|\vp(0)\|^2_{L^2(\Om)^N} \\
\leq C \left( \iint\limits_Q e^{-3s\beta^*}|g|^2 \dxt 
+ \sum_{j=1,j\neq i}^N\int\limits_0^T\int\limits_{\om}e^{-2s\widehat{\beta}-3s\beta^*}\widehat \gamma^7|\vp_j|^2 \dxt \right).
\end{multline}
\end{Lemma}

\textbf{Proof:} We start by an a priori estimate for the Stokes system \eqref{eq:adj-Stokes}. To do this, we introduce a function $\nu\in C^1([0,T])$ such that
$$\nu \equiv 1 \mbox{ in }[0,T/2],\, \nu \equiv 0 \mbox{ in } [3T/4,T].$$
We easily see that $(\nu\vp,\nu\pi)$ satisfies
\begin{equation*}
\left\lbrace \begin{array}{ll}
    -(\nu\vp)_t - \Delta (\nu\vp) + \nabla (\nu\vp) = \nu g - \nu'\vp& \mbox{ in }Q, \\
    \nabla\cdot (\nu\vp) = 0 & \mbox{ in }Q, \\
    (\nu\vp) = 0 & \mbox{ on }\Sigma, \\
    (\nu\vp)(T) = 0 & \mbox{ in }\Om,
 \end{array}\right.
\end{equation*}
thus we have the energy estimate
\begin{equation*}
\|\nu\vp\|^2_{L^2(0,T;H^1(\Om)^N)} + \|\nu\vp\|^2_{L^{\infty}(0,T;L^2(\Om)^N)}\leq C( \|\nu g\|^2_{L^2(Q)^N} + \|\nu' \vp\|^2_{L^2(Q)^N}),
\end{equation*}
from which we readily obtain
\begin{equation*}
\|\vp  \|^2_{L^2(0,T/2;L^2(\Om)^N)}  + \|\vp(0)\|^2_{L^2(\Om)^N}
\leq C( \| g\|^2_{L^2(0,3T/4;L^2(\Om)^N)} + \|\vp\|^2_{L^2(T/2,3T/4;L^2(\Om)^N)}).
\end{equation*}
From this last inequality, and the fact that
$$e^{-3s\beta^*}\geq C>0, \,\forall t\in [0,3T/4] \mbox{ and }e^{-5s\alpha^*}(\xi^*)^4\geq C >0,\, \forall t\in [T/2,3T/4] $$
we have
\begin{multline}\label{eq:proofCarl2-1}
\int\limits_0^{T/2}\int\limits_{\Om} e^{-5s\beta^*}(\gamma^*)^4 |\vp|^2 \dxt +\|\vp(0)\|^2_{L^2(\Om)^N}\\
\leq C\left( \int\limits_0^{3T/4}\int\limits_{\Om} e^{-3s\beta^*}|g|^2 \dxt + \int\limits_{T/2}^{3T/4}\int\limits_{\Om} e^{-5s\alpha^*}(\xi^*)^4|\vp|^2 \dxt \right).
\end{multline}
Note that, since $\alpha=\beta$ in $\Om\times (T/2,T)$, we have:
\begin{equation*}
\begin{split}
\int\limits_{T/2}^T\int\limits_{\Om} e^{-5s\beta^*}(\gamma^*)^4 |\vp|^2 \dxt &= \int\limits_{T/2}^T\int\limits_{\Om} e^{-5s\alpha^*}(\xi^*)^4 |\vp|^2 \dxt \\ &\leq  C \iint\limits_Q e^{-5s\alpha^*}(\xi^*)^4|\vp|^2\dxt,
\end{split}
\end{equation*}
and by the Carleman inequality of Proposition \ref{prop:Carleman}
\begin{multline*}
\int\limits_{T/2}^{T}\int\limits_{\Om} e^{-5s\beta^*}(\gamma^*)^4 |\vp|^2 \dxt \\
\leq C\left( \iint\limits_Q e^{-3s\alpha^*}|g|^2 \dxt +\sum_{j=1,j\neq i}^N \int\limits_{0}^{T}\int\limits_{\om} e^{-2s\widehat{\alpha}-3s\alpha^*}(\widehat \xi)^7|\vp_j|^2 \dxt \right).
\end{multline*}
Since 
$$e^{-3s\beta^*}, e^{-2s\widehat{\beta}-3s\beta^*}\widehat{\gamma}^7 \geq C>0, \,\forall t\in [0,T/2],$$
we can readily get
\begin{multline*}
\int\limits_{T/2}^{T}\int\limits_{\Om} e^{-5s\beta^*}(\gamma^*)^4 |\vp|^2 \dxt \\
\leq C\left( \iint\limits_Q e^{-3s\beta^*}|g|^2 \dxt + \sum_{j=1,j\neq i}^N\int\limits_{0}^{T}\int\limits_{\om} e^{-2s\widehat{\beta}-3s\beta^*}\widehat{\gamma}^7|\vp_j|^2 \dxt \right),
\end{multline*}
which, together with \eqref{eq:proofCarl2-1}, yields \eqref{eq:Carleman2}.

\vskip 1cm

Now we will prove the null controllability of \eqref{eq:Stokes2}. Actually, we will prove the existence of a solution for this problem in an appropriate weighted space. 

Let us set
   \begin{equation*}\label{operatorL}
   Ly=y_t-\Delta y
   \end{equation*}
and let us introduce the space, for $N=2 \mbox{ or }3$ and $i\in \{1,\dots,N\}$,
\begin{equation*}
\begin{array}{l}
E_N^i=\{\, (y,p,v):
 e^{3/2s\beta^*}\,y,\,e^{s\widehat{\beta}+3/2s\beta^*}\widehat \gamma^{-7/2}\,v\1_{\om}\in L^2(Q)^N,\,v_i\equiv 0, \\
 \noalign{\medskip}\hskip1cm e^{3/2 s\beta^*}(\gamma^*)^{-9/8}y \in L^2(0,T;H^2(\Om)^N)\cap L^{\infty}(0,T;V),\\
\noalign{\medskip}\hskip1cm \, e^{5/2s\beta^*}(\gamma^*)^{-2}(Ly+\nabla
p-v\1_{\om}) \in L^2(Q)^N\,\}.
\end{array}
\end{equation*}

It is clear that $E_N^i$ is a Banach space for the following norm:
$$
\begin{array}{l}
\displaystyle \|(y,p,v)\|_{E_N^i}=\left( \|e^{3/2s\beta^*}\,y\|^2_{L^2(Q)^N}
+\|e^{s\widehat{\beta}+3/2s\beta^*}\widehat\gamma^{-7/2}\,v\1_{\om}\|^2_{L^2(Q)^N}\right.\\
 + \|e^{3/2 s\beta^*}(\gamma^*)^{-9/8}\,y\|^2_{L^2(0,T;H^2(\Om)^N)} 
\displaystyle +\|e^{3/2 s\beta^*}(\gamma^*)^{-9/8}y\|^2_{L^\infty(0,T;V)}\\
\left.+\|e^{5/2s\beta^*}(\gamma^*)^{-2}(Ly+\nabla p -v\1_{\om})\|^2_{L^2(Q)^N} \right)^{1/2}
\end{array}
$$

\begin{Remark}
Observe in particular that $(y,p,v)\in E^i_N$ implies $y(T)=0$ in $\Om$. Moreover, the functions belonging to this space posses the interesting following property:
$$e^{5/2s\beta^*}(\gamma^*)^{-2}(y\cdot \nabla)y\in L^2(Q)^N.$$
\end{Remark}

\begin{Proposition}\label{prop:null}
Let $i\in \{1,\dots,N\}$. Assume that 
\begin{equation*}
y^0\in V \mbox{ and }e^{5/2s\beta^*}(\gamma^*)^{-2}f\in L^2(Q)^N.
\end{equation*}
   Then, we can find a control $v$ such that the associated solution $(y,p)$
to \eqref{eq:Stokes2} satisfies $(y,p,v)\in E_N^i$.
   In particular, $v_i\equiv 0$ and $y(T)=0$.
\end{Proposition}

\noindent {\bf Sketch of the proof:}
   The proof of this proposition is very similar to the one of Proposition~2
in~\cite{E&S&O&P} and Proposition 1 in \cite{E&S&O&P-N-1}, so we will just give the main ideas.

Following the arguments in~\cite{FurIma} and~\cite{OlegN-S}, we introduce the space
$$P_0=\{\,(\chi,\sigma)\in C^2(\overline Q)^{N+1}:\nabla\cdot \chi=0,\ \chi=0\ \hbox{on}\
\Sigma\,\}$$
and we consider the following variational problem:
\begin{equation}\label{48p}
a((\widehat \chi,\widehat \sigma),(\chi,\sigma))=\langle G,(\chi,\sigma)\rangle \quad
\forall (\chi,\sigma) \in P_0,
\end{equation}
where we have used the notations
$$\begin{array}{l}
\displaystyle a((\widehat \chi,\widehat \sigma),(\chi,\sigma))=\iint\limits_Q
e^{-3s\beta^*}\,(L^*\widehat \chi+\nabla \widehat \sigma)\cdot(L^*\chi+\nabla \sigma)\,dx\,dt
\\ \noalign{\medskip}\displaystyle
\qquad+\sum_{j=1,j\neq i}^N\int\limits_0^T\int\limits_{\om}
e^{-2s\widehat{\beta}-3s\beta^*}\widehat \gamma^7\,\widehat \chi_j\,\chi_j\,dx\,dt,
\end{array}$$
$$\langle G,(\chi,\sigma)\rangle =\iint\limits_Q  f\cdot \chi
\,\dxt+ \int\limits_{\Omega}y^0\cdot \chi(0)\,dx$$
and $L^*$ is the adjoint operator of $L$, i.e.
   $$
L^*\chi = -\chi_t - \Delta \chi.
   $$

   It is clear that $a(\cdot\,,\cdot):P_0\times P_0\mapsto\R$ is a symmetric,
definite positive bilinear form on $P_0$.
   We denote by $P$ the completion of $P_0$ for the norm induced by
$a(\cdot\,,\cdot)$.
   Then $a(\cdot\,,\cdot)$ is well-defined, continuous and again definite
positive on $P$.
   Furthermore, in view of the Carleman estimate \eqref{eq:Carleman2}, the linear form
$(\chi,\sigma) \mapsto \langle G,(\chi,\sigma)\rangle$ is well-defined and continuous on $P$.
   Hence, from Lax-Milgram's lemma, we deduce that the
variational problem
   \begin{equation}\label{59}
   \left\{
   \begin{array}{l}
   \displaystyle a((\widehat \chi,\widehat \sigma),(\chi,\sigma))=\langle
   G,(\chi,\sigma)\rangle
   \\ \noalign{\medskip}\displaystyle
   \forall (\chi,\sigma) \in P, \quad (\widehat \chi,\widehat \sigma) \in P,
   \end{array}
   \right.
   \end{equation}
   possesses exactly one solution $(\widehat \chi,\widehat \sigma)$.

Let $\widehat y$ and $\widehat v$ be given by 
   \begin{equation*}
   \left\{\begin{array}{ll} \displaystyle \widehat
   y=e^{-3s\beta^*}(L^*\widehat \chi+\nabla \widehat \sigma),&\mbox{ in }Q,
   \\ \noalign{\smallskip}\displaystyle
   \widehat v_j=-e^{-2s\widehat{\beta}-3s\beta^*}\widehat\gamma^7\,\widehat \chi_j\quad (j\neq i),\quad
   \widehat v_i\equiv0&\mbox{ in
   }\om\times (0,T).
   \end{array}\right.
   \end{equation*}
   Then, it is readily seen that they satisfy
 $$
\iint\limits_{Q}e^{3s\beta^*} |\widehat y|^2dxdt
+\sum_{j=1,j\neq i}^N\int\limits_0^T\int\limits_{\om} e^{2s\widehat{\beta}+3s\beta^*}\widehat \gamma^{-7} |\widehat
v_j|^2dxdt = a((\widehat \chi,\widehat \sigma),(\widehat \chi,\widehat \sigma))<+\infty
   $$
and also that $\widehat y$ is, together with some pressure $\widehat p$,
the weak solution (belonging to $L^2(0,T;V)\cap L^{\infty}(0,T;H)$) of the
Stokes system \eqref{eq:Stokes2} for $v=\widehat v$.

It only remains to check that $$e^{3/2s\beta^*}(\gamma^*)^{-9/8}\widehat y\in L^2(0,T;H^2(\Om)^N)\cap
L^{\infty}(0,T;V).$$
  To this end, we define the functions
$$y^*=e^{3/2 s\beta^*}(\gamma^*)^{-9/8}\,\widehat y,
\,p^*=e^{3/2 s\beta^*}(\gamma^*)^{-9/8}\,\widehat p$$ and
$$f^*=e^{3/2 s\beta^*}(\gamma^*)^{-9/8}(f+\widehat v\1_{ \om}).$$ Then
$(y^*,p^*)$ satisfies
   \begin{equation}\label{ystar}
   \left\{\begin{array}{ll} \displaystyle Ly^*+\nabla
   p^*=f^*+(e^{3/2 s\beta^*}(\gamma^*)^{-9/8})_t\,\widehat y&\mbox{ in }Q,\\
   \nabla\cdot y^*=0&\mbox{ in }Q,\\
    \noalign{\smallskip}\displaystyle
   y^*=0&\mbox{ on }\Sigma,
   \\ \noalign{\smallskip}\displaystyle
   y^*(0)=e^{3/2s\beta^*(0)}(\gamma^*(0))^{-9/8}y^0&\mbox{ in }\Omega.
   \end{array}\right.
   \end{equation}
   From the fact that $f^*+(e^{3/2 s\beta^*}(\gamma^*)^{-9/8})_t\,\widehat y \in L^2(Q)^N$ and $y^0\in V$, we
have indeed
   $$
y^*\in L^2(0,T;H^2(\Om)^N)\cap L^{\infty}(0,T;V)
   $$
(see \eqref{eq:regularity1}). This ends the sketch of the proof of Proposition \ref{prop:null}.

\section{Proof of Theorem \ref{teo:nullcontrol}}

In this section we give the proof of Theorem \ref{teo:nullcontrol} using similar arguments to those in \cite{OlegN-S} (see also \cite{E&S&O&P} and \cite{E&S&O&P-N-1}). The result of null controllability for the linear system \eqref{eq:Stokes2} given by Proposition \ref{prop:null} will allow us to apply an inverse mapping theorem. Namely, we will use the following theorem (see \cite{ATF}).

\begin{Theorem}\label{teo:invmap}
Let $B_1$ and $B_2$ be two Banach spaces and let $\mathcal{A}:B_1 \to B_2$ satisfy $\mathcal{A}\in C^1(B_1;B_2)$. Assume that $b_1\in B_1$, $\mathcal{A}(b_1)=b_2$ and that $\mathcal{A}'(b_1):B_1 \to B_2$ is surjective. Then, there exists $\delta >0$ such that, for every $b'\in B_2$ satisfying $\|b'-b_2\|_{B_2}< \delta$, there exists a solution of the equation
$$\mathcal{A} (b) = b',\quad b\in B_1.$$
\end{Theorem}

We apply this theorem setting, for some given $i\in \{1,\dots,N\}$,
 $$B_1 = E_N^i,$$ $$B_2 = L^2(e^{5/2 s\beta^*}(\gamma^*)^{-2}(0,T);L^2(\Om)^N) \times V$$ 
and the operator
$$\mathcal{A}(y,p,v) = (Ly + (y\cdot \nabla)y +\nabla p - v\1_{\om},y(0)) $$
for $(y,p,v)\in E_N^i$.

In order to apply Theorem \ref{teo:invmap}, it remains to check that the operator $\mathcal{A}$ is of class $C^1(B_1;B_2)$. Indeed, notice that all the terms in $\mathcal{A}$ are linear, except for $(y\cdot \nabla)y$. We will prove that the bilinear operator
$$((y_1,p_1,v_1),(y_2,p_2,v_2))\to(y_1\cdot \nabla)y_2$$
is continuous from $B_1\times B_1$ to $ L^2(e^{5/2 s\beta^*}(\gamma^*)^{-2}(0,T);L^2(\Om)^N)$. To do this, notice that $e^{3/2 s\beta^*}(\gamma^*)^{-9/8}y \in L^2(0,T;H^2(\Om)^N)\cap L^{\infty}(0,T;V)$ for any $(y,p,v)\in B_1$, so we have
$$e^{3/2 s\beta^*}(\gamma^*)^{-9/8}y \in L^2(0,T;L^{\infty}(\Om)^N)$$
and  
$$\nabla (e^{3/2 s\beta^*}(\gamma^*)^{-9/8}y) \in L^{\infty}(0,T;L^2(\Om)^N).$$
Consequently, we obtain
\begin{equation*}
\begin{split}
&\|e^{5/2 s\beta^*}(\gamma^*)^{-2}(y_1\cdot \nabla)y_2\|_{L^2(Q)^N} \\ 
&\leq C  \|(e^{3/2 s\beta^*}(\gamma^*)^{-9/8}\,y_1\cdot \nabla)e^{3/2 s\beta^*}(\gamma^*)^{-9/8}\,y_2\|_{L^2(Q)^N} \\
&\leq C \|e^{3/2 s\beta^*}(\gamma^*)^{-9/8}y_1\|_{L^2(0,T;L^{\infty}(\Om)^N)}\, \|e^{3/2 s\beta^*}(\gamma^*)^{-9/8}y_2\|_{L^{\infty}(0,T;H^1(\Om)^N)}.
\end{split} 
\end{equation*}

Notice that $\mathcal{A}'(0,0,0):B_1\to B_2$ is given by
$$\mathcal{A}'(0,0,0)(y,p,v) = (Ly + \nabla p,y(0)),\, \forall (y,p,v)\in B_1,$$
so this functional is surjective in view of the null controllability result for the linear system \eqref{eq:Stokes2} given by Proposition \ref{prop:null}.

We are now able to apply Theorem \ref{teo:invmap} for $b_1=(0,0,0)$ and $b_2=(0,0)$. In particular, this gives  the existence of a positive number $\delta$ such that, if $\|y(0)\|_V\leq \delta$, then we can find a control $v$ satisfying $v_i\equiv 0$, for some given $i\in \{1,\dots,N\}$, such that the associated solution $(y,p)$ to \eqref{eq:NS} satisfies $y(T)=0$ in $\Om$.

This concludes the proof of Theorem \ref{teo:nullcontrol}.


\end{document}